\documentclass[12pt,a4paper,leqno]{article}

\usepackage[T1]{fontenc}
\usepackage[francais]{babel}
\usepackage{amsmath}
\usepackage{amssymb}
\usepackage{url}

\def\iy{\infty}

\def\al{\alpha}

\def\ga{\gamma}

\def\de{\delta}
\def\De{\Delta}

\def\th{\theta}
\def\Th{\Theta}

\def\vep{\varepsilon}

\def\om{\omega}
\def\Om{\Omega}
\def\si{\sigma}

\def\ca{\mathcal{A}}

\def\cd{\mathcal{D}}

\def\cg{\mathcal{G}}

\def\cn{\mathcal{N}}
\def\co{\mathcal{O}}
\def\cp{\mathcal{P}}

\def\cR{\mathcal{R}}
\def\cs{\mathcal{S}}

\def\Llr{\Longleftrightarrow}

\def\mk{\medskip}

\def\Llr{\Longleftrightarrow}

\newtheorem{theorem}{Theorem}
\newenvironment{dem}{\smallskip\noindent{\it Proof.}%
{\nopagebreak[0]}}%
{\nopagebreak[0]\hfill$\Box$ \medskip}%
\newtheorem{lem}{Lemma}
\newtheorem{coro}{Corollary}

\title{On SA, CA, and GA numbers}
\author{Geoffrey Caveney\\
7455 North Greenview \#426,
Chicago, IL 60626, USA\\
E-mail: rokirovka@gmail.com\\
\and
Jean-Louis Nicolas\\
Universit\'e de Lyon; CNRS; Universit\'e Lyon 1;\\
Institut Camille Jordan, Math\'ematiques,\\
21 Avenue Claude Bernard, F-69622 Villeurbanne cedex, France\\
E-mail: nicolas@math.univ-lyon1.fr\\
\and
Jonathan Sondow\\
209 West 97th Street \#6F,
New York, NY 10025, USA\\
E-mail: jsondow@alumni.princeton.edu\\
Tel.: +1-646-306-1909}

\begin{document}
\date{}
\maketitle

\noindent{\bf Abstract} Gronwall's function $G$ is defined for $n>1$ by
$G(n)=\frac{\sigma(n)}{n \log\log n}$
where $\sigma(n)$ is the sum of the divisors of $n$. We call an integer $N>1$ a~\emph{GA1 number} if $N$~is composite and $G(N) \ge G(N/p)$ for all prime factors~$p$ of~$N$. We say that $N$ is a \emph{GA2 number} if $G(N) \ge G(aN)$ for all multiples $aN$ of~$N$. In arXiv 1110.5078, we used Robin's and Gronwall's theorems on~$G$ to prove that the Riemann Hypothesis (RH) is true if and only if $4$~is the only number that is both GA1 and GA2. Here, we study GA1 numbers and GA2 numbers separately. We compare them with superabundant (SA) and colossally abundant (CA) numbers (first studied by Ramanujan). We give algorithms for computing GA1 numbers; the smallest one with more than two prime factors is 183783600, while the smallest odd one is 1058462574572984015114271643676625. We find nineteen GA2 numbers $\le5040$, and prove that a GA2 number $N>5040$ exists if and only if RH is false, in which case $N$ is even and $>10^{8576}$.

\mk

\mk

\section{Introduction}\label{secIntr}

The \emph{sum-of-divisors function} $\sigma$ is defined by
$$\sigma(n):=\sum_{d \mid n}d.$$
For example, $\sigma(4)=7$.

In 1913, Gronwall \cite{gronwall} found the maximal order of~$\sigma$.

\begin{theorem}[Gronwall] \label{THM: gronwall}
The function
$$G(n):=\frac{\sigma(n)}{n \log\log n} \qquad(n>1)$$
satisfies
$$\limsup_{n\to\infty} G(n) =  e^\gamma = 1.78107\dotso,$$
where $\gamma$ is the Euler-Mascheroni constant.
\end{theorem}

In 1915, Ramanujan proved an asymptotic inequality for Gronwall's
function $G$, assuming the Riemann Hypothesis (RH). Ramanujan's result
was shown in the second part of his thesis. The first part was
published in 1915 \cite{ramanujan15} while the second part was 
not published until much later, in 1997 \cite{ramanujan97}.

\begin{theorem}[Ramanujan] \label{THM: ramanujan}
If the Riemann Hypothesis is true, then
\begin{equation*}
G(n) < e^\gamma \qquad (n \gg 1).
\end{equation*}
\end{theorem}
Here, $n\gg1$ means for all sufficiently large $n$.

In 1984, without being aware of Ramanujan's theorem, Robin \cite{robin} proved that a stronger statement about the function~$G$ is \emph{equivalent} to RH.

\begin{theorem}[Robin] \label{THM: robin}
The Riemann Hypothesis is true if and only if
\begin{equation}
G(n) < e^\gamma \qquad (n > 5040). \label{EQ: robin}
\end{equation}
\end{theorem}

The condition~\eqref{EQ: robin} is called \emph{Robin's
  inequality}. Table~\ref{TABLE: S and G(m)} gives the twenty-six known
numbers~$r$ for which the reverse inequality $G(r) \ge  e^\gamma$
holds (see \cite[Sequence A067698]{oeis}), together with the value of $G(r)$ (truncated). (The ``$a(r)$'' column is explained in \S\ref{secGA2}, and the ``$Q(r)$'' column in \S\ref{subcompG}.)

\begin{table}[htdp]
\begin{center}
\begin{tabular}{r@{     }c@{     }c@{     }c@{     }c@{     }l@{     }c@{     }c@{     }c@{     }r}
\hline
$r\ $ & SA & CA & GA1  & GA2 &  Factorization &  $\sigma(r)/r$ & $G(r)$ & $a(r)$ & $Q(r)$\\
\hline
$3$  &  & & & \checkmark &   3 & 1.333 & $14.177\ \ $ &0 \\
$4$ & \checkmark & & \checkmark & \checkmark&  $2^2$  & 1.750 & 5.357 &0 &$-0.763$ \\
$5$   &  & & & \checkmark & 5  &  1.200 & 2.521 &0 \\
$6$   & \checkmark &\checkmark & & \checkmark &   $2\cdot3$ & 2.000 & 3.429 &0&$4.134$ \\
$8$ &   &  &  &\checkmark & $2^3$ &   1.875 & 2.561 &0&$2.091$  \\
$9$  & & & & & $3^2$ &  1.444 & 1.834 & 4&$7.726$\\
$10$  &  &  && \checkmark &  $2\cdot5$  & 1.800 & 2.158  &0&$1.168$\\
$12$ & \checkmark &\checkmark  & & \checkmark &  $2^2\cdot3$  & 2.333 & 2.563  &0&$2.090$\\
$16$  & & & & &  $2^4$ & 1.937 & 1.899 & 3&$1.348$ \\
$18$  &  &  &   &\checkmark & $2\cdot3^2$   & 2.166 & 2.041&0&$1.679$ \\
$20$  &  &  & &  &$2^2\cdot5$ &  2.100 & 1.913 & 3 &$2.799$\\
$24$   & \checkmark & &  &\checkmark &   $2^3\cdot3$ & 2.500 & 2.162  &0&$1.185$ \\
$30$ &   & &&   & $2\cdot3\cdot5$   & 2.400 & 1.960 & 2 &$1.749$\\
$36$  &  \checkmark&  &&\checkmark &  $2^2\cdot3^2$ & 2.527 & 1.980  &0&$1.294$\\
$48$ &   \checkmark & & &\checkmark &   $2^4\cdot3$  & 2.583 & 1.908&0&$1.132$ \\
$60$ &  \checkmark  &\checkmark  & & \checkmark &   $2^2\cdot3\cdot5$  & 2.800 & 1.986 &0&$1.290$\\
$72$ &   & & & \checkmark & $2^3\cdot3^2$   & 2.708 & 1.863 &0&$1.160$  \\
$84$  &  &  & & & $2^2\cdot3\cdot7$  & 2.666 & 1.791 & 10 &$1.430$ \\
$120$ & \checkmark &\checkmark  & & \checkmark &  $2^3\cdot3\cdot5$   & 3.000 & 1.915 &0&$1.128$\\
$180$  & \checkmark & &  &\checkmark  &  $2^2\cdot3^2\cdot5$  & 3.033 & 1.841&0&$1.078$\\
$240$  & \checkmark & & & \checkmark &   $2^4\cdot3\cdot5$ & 3.100 & 1.822  &0&$1.051$\\
$360$   & \checkmark &\checkmark & & \checkmark  &  $2^3\cdot3^2\cdot5$ & 3.250 & 1.833 &0&$1.044$\\
$720$  &\checkmark && & &   $2^4\cdot3^2\cdot5$  & 3.358 & 1.782 & 7&$1.028$\\
$840$ & \checkmark & &  &&   $2^3\cdot3\cdot5\cdot7$& 3.428 & 1.797 & 3 &$1.065$ \\
$2520$  & \checkmark &\checkmark & & \checkmark &   $2^3\cdot3^2\cdot5\cdot7$  & 3.714 & 1.804  &0&$1.015$\\
$5040$  &\checkmark & \checkmark &  &\checkmark &   $2^4\cdot3^2\cdot5\cdot7$& 3.838 & 1.790&0&$1.007$\\
\hline
\end{tabular}
\end{center}
\caption{The set $\cR=\{r\le5040:G(r) \ge  e^\gamma= 1.781\dotso\}$, which contains the subset $\cR_2=\{N\le5040:N\text{ is GA}2\}.$}
\label{TABLE: S and G(m)}
\end{table}

In \cite{robin} Robin also proved, unconditionally, that
\begin{equation}
G(n) \le e^\gamma + \frac{0.6482\dotso}{(\log\log n)^2} \qquad (n > 1)\label{EQ: bound}
\end{equation}
with equality for $n=12$. This refines the inequality $\limsup_{n\to\infty} G(n) \le  e^\gamma$ from Gronwall's theorem.

Recently, the authors \cite{CNS} used Robin's results to derive another reformulation of RH. Before recalling its statement, we give three definitions and an example.

\mk
A positive integer $N$ is a \emph{GA1 number} if $N$ is composite and the inequality
$$G(N) \ge G(N/p)$$
holds for all prime factors $p$ of $N$. The first few GA1 numbers are
$$N=4, 14, 22, 26, 34, 38, 46, 58, 62, 74, 82, 86, 94, 106, 118, 122, 134, 142,\dotsc$$
 (see \cite[Sequence A197638]{oeis}), and (see \S\ref{secodd}) {\it the smallest odd GA1 number is}
$$N=1058462574572984015114271643676625.$$

\mk
An integer $N>1$ is a \emph{GA2 number} if$$G(N) \ge G(aN)$$for all multiples $aN$ of $N$. The nineteen known GA2 numbers (see Theorem~\ref{thmGA2} and \cite[Sequence A197369]{oeis}) are
$$N=3,4,5,6,8,10,12,18,24,36,48,60,72,120,180,240,360,2520,5040.$$
\emph{Every GA2 number $>5$ is even.} (\emph{Proof.} 
If $N$ is odd, then $\sigma(2N)=3\mspace{1mu}\sigma(N)$,
and if $N$ is also GA2, we get
$$\frac32 \le \frac{3G(N)}{2G(2N)} = \frac{\log\log 2N}{\log\log N}$$
which implies $N<7$.)

\mk
Finally, a composite number is \emph{extraordinary} if it is both GA1 and GA2.

\mk
For example, the smallest extraordinary number is $4.$ To see this, we first compute $G(4)=5.357\dotso.$ Then, as $G(2)<0$, it follows that $4$ is a GA1 number. Since Robin's unconditional bound \eqref{EQ: bound} implies
\begin{equation*}
G(n) < e^\gamma + \frac{0.6483}{(\log\log 5)^2} = 4.643\dotso < G(4) \qquad(n\ge5),
\end{equation*}
we get that $4$ is also GA2. Thus $4$ is an extraordinary number.

We can now recall our results from \cite[Theorem~6 and Corollary~8]{CNS}.

\begin{theorem}[Caveney-Nicolas-Sondow] \label{THM: Caveney}
\emph{(i)}. The Riemann Hypothesis is true if and only if $4$ is the only extraordinary number.

\noindent\emph{(ii)}. If there is any counterexample to Robin's inequality, then the maximum $\mu:=\max\{G(n):n > 5040\}$ exists and the least number $N>5040$ with $G(N)=\mu$ is extraordinary.
\end{theorem}

If there exists an extraordinary number $N>4$, then $N$ is even (as $5$~is not GA1, and no GA2 number $>5$ is odd) and $N>10^{8576}$ (since  no GA1 number lies in the interval $[5,5040]$, and no GA2 number lies in $[5041,10^{8576}]$---see Corollary~\ref{CORlowerbound}).

In the present paper, we study GA1 numbers and GA2 numbers separately.

Preliminary facts about GA1 numbers and GA2 numbers were given in  \cite{CNS}. We recall two of them and make a definition.

\mk
\noindent{\bf Fact 1} (proved by elementary methods in \cite[\S5]{CNS}). \emph{The GA1 numbers with exactly two $($not necessarily distinct$)$ prime factors are precisely $4$ and $2p$, for primes} $p\ge7$.

\mk
We call such GA1 numbers {\it improper}, while GA1 numbers with at least three (not
necessarily distinct) prime factors will be called {\it proper}.

\emph{The smallest proper GA1 number is} $\nu:=183783600$ (see \S\ref{secodd} and \cite[Sequence A201557]{oeis}). The number $\nu$ was mentioned in \cite[equation (3)]{CNS} as an example of a (proper) GA1 number that is not a GA2 number (because $G(\nu)<G(19\nu)$).

\mk
\noindent{\bf Fact 2} (see \cite[Lemma 10]{CNS}). \emph{If $n_0$ is a positive integer,
then
$$\limsup_{a\to \iy} G(an_0)=e^\ga,$$
which yields the implication}
\begin{equation}\label{GA2>ega}
N\ \text{\it is GA2} \implies G(N) \geq e^\ga.
\end{equation}

An application is an alternate proof that any GA2 number $N>5$ is even. Namely, as $7$ and $9$ are not GA2, and as Theorem~2 in \cite{clms} says that an integer $n>9$ is even if $G(n) \geq e^\ga$, the result follows from \eqref{GA2>ega}.

By the method of \cite[\S5]{CNS}, one can prove two additional properties of GA1 numbers.

\mk
\noindent{\bf Fact 3}. \emph{The only prime power GA1 number $N=p^k$ is} $N=4.$

\mk
\noindent{\bf Fact 4}. \emph{A product of three distinct primes $p_1p_2p_3$ cannot be a GA1 number.} (See \S\ref{secOm(N)} for a more general result proved by other methods.)

\mk
The rest of the paper is organized as follows. The next subsection establishes notation. In \S\ref{secCASA} we recall the definitions of superabundant (SA) and colossally abundant (CA) numbers and review some of their properties. In \S\ref{seclem} we prove six lemmas needed later. In \S\ref{secGA2} we give an analog of Theorem~\ref{THM: Caveney} for GA2 numbers; in particular, \emph{if RH is false, then infinitely many GA2 numbers exist, and any number~$N> 5040$ for which $G(N)=\max\{G(n):n > 5040\}$ is both GA2 and CA.} In the final four sections we study proper GA1 numbers: \S\ref{secCAGA1} compares them with SA and CA numbers, \S\ref{secpr} is concerned with their prime factors, \S\ref{seccomp} gives algorithms for computing them, and \S\ref{secconcl} estimates the number of them up to $x$.

\subsection{Notation} \label{SEC: notation}

We let $p$ always denote a prime.

Let $v_p(n)$ denote \emph{the exponent on $p$ in the prime
factorization}
$$n=\prod_{p} p^{v_p(n)}.$$

For $n\geq 1$, we denote \emph{the number of prime factors of
$n$ counted with multiplicity} by
$$\Om(n):=\sum_{p} v_p(n).$$

For $n>1$, we denote \emph{the largest prime factor of} $n$ by
$$P(n):=\max\{p:p\mid n\} = \max\{p:v_p(n)>0\}.$$

As usual, Chebychev's function $\th$ is defined as
$$\theta(x):=\sum_{p\,\le\, x} \log p.$$

\section{Review of properties of SA and CA numbers}\label{secCASA}

 Superabundant and colossally abundant
 numbers were first introduced by Ramanujan, who called them
 generalized  highly composite and generalized super highly composite numbers, respectively (cf. \cite[\S
59]{ramanujan97}). They were rediscovered later by Alaoglu and Erd\H
{o}s \cite{ae}.

A \emph{superabundant} (\emph{SA}) number is a positive integer $N$ such that
\begin{equation*} 
\frac{\sigma(N)}{N} > \frac{\sigma(n)}{n} \qquad (0<n<N).
\end{equation*}
The first few SA numbers are (see \cite[Sequence A004394]{oeis})
$$N=1, 2, 4, 6, 12, 24, 36, 48, 60, 120, 180, 240, 360, 720, 840, 1260, 1680,  \dotso.$$

A \emph{colossally abundant} (\emph{CA}) number is a positive integer $N$ for which there exists an exponent $\vep>0$ such that
\begin{equation}\label{defCA}
\frac{\sigma(N)}{N^{1+\vep}} \ge \frac{\sigma(n)}{n^{1+\vep}} \qquad (n>1).
\end{equation}
Such an exponent $\vep$ is called a \emph{parameter} of $N$. The
sequence of CA numbers
(compare \cite[Sequence A004490]{oeis}) begins
$$N=1,2, 6, 12, 60, 120, 360, 2520, 5040, 55440, 720720, 1441440, 4324320, \dotso.$$
From \eqref{defCA}, it is easy to show that \emph{every CA number is also SA}. 

Now let $N$ denote an SA or CA number. Then (see \cite[Theorems~1 and~3]{ae} or \cite[\S 59]{ramanujan97})
\begin{equation}\label{k}
N = 2^{k_2}\cdot3^{k_3}\cdot5^{k_5}\dotsb p^{k_p}  \quad  \implies
\quad k_2 \ge k_3 \ge k_5 \ge \dotsb \ge k_p
\end{equation}
with $k_p=1$ unless $N=4$ or $36$, and \cite[Theorem 7]{ae}
\begin{equation}\label{PNeqlogN}
p=P(N) \sim \log N \quad (N\to\infty).
\end{equation}

We recall some properties of CA numbers (see \cite{ae, briggs,
  CNS, EN, lagarias, ramanujan97, robin, robinsem}).

Note first that for any fixed positive integer $k$, the quantity
$$F(t,k):=\frac{\log \left(1+\dfrac{1}{t+t^2+\dotsb +t^k}\right)}{\log t}$$
is decreasing on the interval $1<t<\infty$, and the function $t\mapsto F(t,k)$ maps the interval onto the positive real numbers. Hence, given $\vep > 0$, we may define $x_k = x_k(\vep) > 1$ by
\begin{equation}\label{defxk}
F(x_k,k)=\vep.
\end{equation} 
(See \cite[p. 189]{robin} and \cite[\S 61 and \S 69] {ramanujan97}.)
In particular, when $k=1$ we set $x=x_1=x_1(\vep)$, so that
\begin{equation}\label{defx}
F(x,1)=F(x_1,1)=\frac{\log \left(1+\dfrac{1}{x}\right)}{\log x}=\vep.
\end{equation}

It is convenient to set $x_0=+\iy$. From the decreasingness of
$F(t,k)$ with respect to both $t$ and $k$, it follows that the sequence
$(x_k)_{k\geq 0}$ is decreasing.

If $N$ is a CA number of parameter $\vep$ and $p$ divides $N$ with
$v_p(N)=k$, then applying \eqref{defCA} with $n=Np$ yields
\begin{equation*} 
\vep \geq F(p,k+1)\qquad \text{ i.e. } \quad p\geq x_{k+1}
\end{equation*} 
while, if $k > 0$, applying \eqref{defCA} with $n=N/p$ yields
\begin{equation}\label{eps<}
\vep \leq F(p,k)\qquad \text{ i.e. } \quad p\leq x_{k}.
\end{equation} 
Let $K$ be the largest integer such that $x_K\geq 2$. Then from
\eqref{eps<}, for all $p$'s we have $2\leq p \leq x_k$ and 
\begin{equation*} 
k=v_p(N) \leq K.
\end{equation*} 

Now define the set
$${\cal E}:=\{F(p,k): p\text{ is prime and }k \ge1\}.$$
Its largest element is
$$\max{\cal E} = F(2,1) = \frac{\log(3/2)}{\log2} = 0.5849\dotso,$$
and its infimum is
$$\inf{\cal E} = \lim_{k\to\infty} F(p,k) = 0$$
for any fixed prime $p$. 

If $\vep \notin {\cal E}$, then no $x_k$ is a prime number and there exists
a unique CA number $N=N(\vep)$ of parameter $\vep;$ moreover, $N$ is given by either of the equivalent formulas
$$N=\prod_{p < x}p^{k_p} \qquad \text{ with } x_{k_p+1} < p < x_{k_p}$$
or
\begin{equation}\label{N<}
N= \prod_{k=1}^K \;\; \prod_{p < x_k}  p.
\end{equation}
In particular, if $\vep>\max{\cal E}$, then $x=x_1 <2$, $K=0$ and $N(\vep)=1$.

If $\vep \in {\cal E}$, then some $x_k$ is prime, and it is
highly probable that only one $x_k$ is prime.
But (see \cite[Proposition 4]{EN}), from the theorem of six
exponentials it is only possible to show that at most two $x_k$'s are
prime. (Compare \cite[p.~538]{lagarias}.) Therefore there are either two or four CA numbers of parameter
$\vep,$ defined by
\begin{equation}\label{N=}
N= \prod_{k=1}^K \;\;\prod_{\substack{p \leq x_k\\
\text{or }\\p < x_k}}  p.
\end{equation}
Here, if $ x_k$ is a prime $p$ for some $k$, then $p$ may or may not be a factor in the inner
product. (This can occur for at most two values of $k$.) In other words, if $x_{k-1} < p < x_k$, then the exponent $v_p(N)$ of $p$ in $N$ is~$k$, while if $p = x_k$, the exponent may be $k$ or $k-1$.
In particular, if $N$ is the \emph{largest} CA number of parameter $\vep$, then
\begin{equation}\label{p||N}
F(p,1)=\vep\quad\implies \quad P(N)=p.
\end{equation}
Note that, since if $\vep \notin {\cal E}$, then $x_k$ is not prime,
formula \eqref{N<} gives the same value as \eqref{N=}. Therefore,
for any $\vep$, formula \eqref{N=} gives all the possible values of a
CA number $N$ of parameter $\vep$. 
(Thus $N$ is a product of ``primorials'' \cite[Sequence A002110]{oeis}.)

\section{Six lemmas}\label{seclem}

The case $k=2$ of the following lemma was proved in \cite[p. 190]{robin}.

\begin{lem}\label{lemzx}
For $k\ge 2$, we have the upper bound
$$x_k < (kx)^{1/k}.$$
\end{lem}

\begin{dem}
Since the function $t\mapsto F(t,k)$ is strictly decreasing on
$1<t<\infty$, to prove $x_k < z:=(kx)^{1/k}$, it suffices to show
$F(z,k)<F(x_k,k)$. As \eqref{defxk} and \eqref{defx} imply
$F(x_k,k)=\vep=F(x,1)$, this reduces to showing $F(z,k)<F(x,1)$.

Since $z>1$ and $k\ge2$, we have
\begin{eqnarray*}
F(z,k)&=&\log \left(1+\frac{1}{z+z^2+\dotsb +z^k}\right)\frac{1}{\log z}\\
& < &\frac{1}{(z+z^2+\dotsb +z^k)\log z}=\frac{k}{(z+z^2+\dotsb +z^k)\log kx}\\
& < &\frac{k}{(k-1+z^k)\log x} \le \frac{k}{\left(\frac k2+kx\right)\log x}=\frac{1}{\left(x+\frac12\right)\log x}\\
&<& \log \left(1+\frac{1}{x}\right)\frac{1}{\log x} = F(x,1),
\end{eqnarray*}
using the lower bound
$\log \left(1+\frac{1}{t}\right) > \left(t+\frac12\right)^{-1}$, valid for $t>0$. This proves the desired inequality.
\end{dem}

\mk

In the proof of Theorem \ref{thmGA2} (iii), we will  need the 
following result (see \cite[Lemma 4]{NR}).

\begin{lem}\label{lemNR}
Given a CA number $N_0$ of parameter $\vep_0$, let $N > N_0$ be a
number satisfying 
\begin{equation}\label{lemNR1}
n\geq N_0 \quad \implies \quad \frac{\si(n)}{n^{1+\vep}} \leq \frac{\si(N)}{N^{1+\vep}}
\end{equation}
for some fixed $\vep> 0$. Then $N$ is CA of parameter $\vep$.
\end{lem}

\begin{dem}
Since $N_0$ is CA of parameter $\vep_0$, we have
$$\frac{\si(N)}{\si(N_0)} \leq \left( \frac{N}{N_0}\right)^{1+\vep_0}.$$
On the other hand, \eqref{lemNR1} yields
$$\frac{\si(N)}{\si(N_0)} \geq \left( \frac{N}{N_0}\right)^{1+\vep}.$$
Hence $\vep \leq \vep_0$.

In view of \eqref{lemNR1}, to prove that $N$ is CA of
parameter $\vep$, we only need to show that
$$n < N_0 \quad \implies \quad \frac{\si(n)}{n^{1+\vep}} \leq
\frac{\si(N)}{N^{1+\vep}} \cdot$$
If $n < N_0$, then since $N_0$ is CA and
\eqref{lemNR1} holds, we have
$$\frac{\si(n)}{n^{1+\vep}} =
\frac{\si(n)n^{\vep_0-\vep}}{n^{1+\vep_0}} \leq
\frac{\si(N_0)n^{\vep_0-\vep}}{N_0^{1+\vep_0}} \leq
\frac{\si(N_0)N_0^{\vep_0-\vep}}{N_0^{1+\vep_0}}
=\frac{\si(N_0)}{N_0^{1+\vep}} \leq \frac{\si(N)}{N^{1+\vep}}.$$
This completes the proof of Lemma \ref{lemNR}.
\end{dem}

\mk

The next lemma provides an estimate for a CA number of parameter $\vep$.

\begin{lem}\label{lemCA}
Let $N$ be a CA number of parameter
$\vep < F(2,1)=\log(3/2)/\log2$
and define $x=x(\vep)$ by \eqref{defx}.\\
\emph{(i).} Then
$$\log N \le \th(x)+c \sqrt x$$
for some constant $c>0$.

\noindent\emph{(ii).} 
Moreover, if $N$ is the largest CA number
\footnote{Note that Ramanujan's definition of CA number of parameter $\vep$ in \cite{ramanujan97}
  is not exactly the same as that of Robin in \cite[pp. 189--190]
  {robin}. Ramanujan's definition corresponds to the
  \emph{largest} CA number of parameter $\vep$ for Robin.}
of parameter $\vep$, then
$$\th(x) \leq \log N \leq \th(x)+c \sqrt x.$$
\end{lem}

\begin{dem}
(i). It follows from formula \eqref{N=} for $N$
that if $x_k$ is defined by \eqref{defxk}, then 
\begin{equation}\label{logN<sum}
\log N \le \th(x_1)+\th(x_2)+\dotsb+\th(x_K),
\end{equation}
where $K$ is the largest integer such that $x_K \geq 2$. (Note that
$v_2(N)=K$ or $K-1$, and that $\vep < F(2,1)$ implies $x > 2$ and $K\geq 1$.) 

As $t\mapsto F(t,k)$ is decreasing and \eqref{defxk} holds, we
have
$$F(2,K) \geq F(x_K,K) =\vep =F(x_{K+1},K+1) > F(2,K+1).$$
On the other hand,
$$F(2,K)=\log\left(1+\frac{1}{2^{K+1}-2}\right)\frac{1}{\log 2}\
<\frac{1}{(2^{K+1}-2) \log 2} \leq \frac{1}{2^{K} \log 2} <\frac{2}{2^K}$$
and, from \eqref{defx},
$$\vep=\frac{\log(1+\frac1x)}{\log x} > \frac{1}{(x+1)\log x} \geq
\frac{1}{(x+1)(x-1)} > \frac{1}{x^2}.$$
Thus
$$\frac{2}{2^K} > F(2,K) \geq \vep > \frac{1}{x^2},$$
implying
\begin{equation}\label{K<logx}
K < 1+\frac{2}{\log 2} \log x.
\end{equation}

Since $k\mapsto x_k$ is decreasing, from \eqref{logN<sum} we have 
(compare \cite[equation (368)]{ramanujan97})
$$\log N \le \th(x_1)+\th(x_2)+K\th(x_3).$$
Using $x_2 \le \sqrt{2x}$ and $x_3 \le \sqrt[3]{3x}$ (from Lemma \ref{lemzx}),
together with \eqref{K<logx} and the Prime Number Theorem in the form $\th(t) \sim t$, we deduce (i).

\mk

\noindent(ii). 
From \eqref{N=}, the largest CA number of parameter $\vep$ is
$$N= \prod_{k=1}^K \;\;\prod_{p \leq x_k}  p$$
which implies $\th(x)\leq \log N$,
and (ii) follows from (i).
\end{dem}

In the next lemma, we recall the oscillations of Chebychev's function
$\th$ studied by Littlewood.

\begin{lem}\label{lemP}
\mbox{}
There exists a constant $c>0$ such that for infinitely many primes~$p$
we have
\begin{equation}\label{P1}
\theta(p) < p-c\sqrt p \log \log \log p,
\end{equation}
and for infinitely many other primes $p$ we have
\begin{equation}\label{P2}
\theta(p) > p+c \sqrt{p}\log \log \log p.
\end{equation}
\end{lem}

\begin{dem}
From Littlewood's theorem (see \cite{Lit}), 
we know that there exists a constant $c'>0$ such that for a sequence of
values of $x$ going to infinity we have
\begin{equation}\label{P3}
\theta(x) < x - c'\sqrt x \log \log \log x,
\end{equation}
and for a sequence of
values of $x'$ going to infinity  we have
\begin{equation}\label{P4}
x' + c'\sqrt{x'} \log\log\log x' < \theta(x').
\end{equation}
Let us suppose first that $x$ is large enough and satisfies \eqref{P3}. If $x=p$
is prime, then \eqref{P3} implies \eqref{P1}. Now assume $x$ is not prime, and let $p$ be the
prime following $x$. As the function
$t \mapsto t - c \sqrt t \log\log\log t$ is increasing, we get
\begin{align*}
\theta(p) = 
\theta(x) + \log p < \ &x - c' \sqrt x \log\log\log x + \log p\\
<  \ &p - c' \sqrt p\log\log\log p + \log p,
\end{align*}
which implies \eqref{P1} with $c < c'$ for $x$ large enough.

The proof of \eqref{P2} is easier. Let $x'$ satisfy 
\eqref{P4} and choose the largest prime
$p\le x'$. For $c \leq c'$, we have
\[
\theta(p) = \theta(x') > x' + c' \sqrt{x'} \log\log\log x' >p+c \sqrt p \log \log \log p,
\]
which proves \eqref{P2}.
\end{dem}

\begin{lem}\label{lemmajth}
Chebychev's function $\th(x)$ satisfies
$$\th(x) \leq (1+\al)x,$$
where
$$\al = \al(x):=
\begin{cases}
\quad 0 & \text{ if } x \leq 8\cdot 10^{11},\\
\dfrac{1}{36260} < 0.000028 & \text{ otherwise.} 
\end{cases}
$$
\end{lem}

\begin{dem}
Schoenfeld (cf. \cite[p. 360]{Sch}) proved $\th(x) \leq
1.000081\,x$ for all $x$, and he mentioned that Brent
had checked that $\th(x) < x$ for $x < 10^{11}$. The
stronger results stated here are due to Dusart---see 
\cite[p. 2 and Table 6.6]{Dus}.
\end{dem}

\begin{lem}\label{lemg}
Let $\vep$ be a positive real number. For $t > e$, let us set
\begin{equation}\label{g}
g(t)=g_\vep(t):=\vep \log t- \log \log \log t.
\end{equation}
Then there exists a unique real number $t_0 =t_0(\vep)> e$ such that
\begin{equation} \label{epst0}
\frac{1}{\log t_0\, \log \log t_0}=\vep.
\end{equation}
Moreover, $g(t)$ is decreasing for $e < t < t_0$ and increasing for $t> t_0$.
\end{lem}

\begin{dem}
The derivative of $g$ is
$$g'(t)=\frac 1t \left( \vep-\frac{1}{\log t\;\log\log t}\right).$$
For $t > e$, both $\log t$ and $\log \log t$ are positive and increasing, and the function $t\mapsto 1/(\log t\; \log\log t)$ is a decreasing bijection from $(e,+\iy)$ onto $(0,+\iy)$. Therefore, one can define $t_0 > e$ by \eqref{epst0}. 

Then we have $g'(t) < 0$ for $e < t < t_0$, and $g'(t) > 0$ for $t > t_0$, which completes the proof of Lemma \ref{lemg}.
\end{dem}

\section{GA2 numbers}\label{secGA2}

We first study GA2 numbers. Compare the following result on them with Theorem~\ref{THM: Caveney} on extraordinary numbers.

\begin{theorem}\label{thmGA2}
\emph{(i)}. The set of GA2 numbers $\le5040$ is
\begin{align*}
\cR_2:=\{3,4,5,6,8,10,12,18,24,36,48,60,72,120,180,240,360,2520,5040\}.
\end{align*}
\noindent\emph{(ii)}. If the Riemann Hypothesis is true, then no GA2 number exceeds $5040$.

\noindent\emph{(iii)}. If the Riemann Hypothesis is false, then infinitely many GA2 numbers exist; moreover, the inequality
$$\mu:=\max\{G(n):n > 5040\} > e^\ga$$ 
holds, and any integer $A> 5040$ for which
$G(A)=\mu$ is both GA2 and CA.
\end{theorem}

\begin{dem}
(i). Setting
$$\cR':=\{N\le5040:N\text{ is GA}2\},$$
we have to prove that $\cR' = \cR_2$.

To show $\cR' \subset \cR_2$, choose $N\in\cR'$. From \eqref{GA2>ega}, we have $G(N) \geq e^{\ga},$ so that
\begin{align*}
N\in\cR:=\{r&\le5040 : G(r) \ge  e^\gamma\}\\
=\{3&, 4, 5, 6, 8, 9, 10, 12, 16, 18, 20, 24, 30, 36, 48, 60, 72, 84, 120, 180, \\
2&40, 360, 720, 840, 2520, 5040\},
\end{align*}
by calculating the ``$r$'' column of Table~\ref{TABLE: S and G(m)}. To show that $N$ belongs to the subset $\cR_2\subset\cR$, define for $r\in\cR$ the integer
$$a(r) := 
\begin{cases}
\min \ca_r, & \text{if }\ca_r:=\{a:G(ar) > G(r),\ ar\in\cR\}\neq\emptyset,\\
0,& \text{if }\ca_r=\emptyset.\\
\end{cases}$$
A computation (see the ``$a(r)$'' column of Table \ref{TABLE: S and G(m)}) shows that
\begin{equation} \label{EQ:R2'}
\{r\in\cR :  \ca_r\neq\emptyset\} = \{9,16,20,30,84,720,840\}.
\end{equation}
Since $N$ is GA2, it must lie in the complement
\begin{equation} \label{EQ:R2}
\cR\setminus \{9,16,20,30,84,720,840\} = \cR_2.
\end{equation}
This shows $\cR' \subset \cR_2$.

To prove $\cR_2 \subset \cR'$, choose $r\in\cR_2$. To get $r\in\cR'$, we need to show that $G(r) \ge G(ar)$, for any multiple $ar$ of $r$. We consider two cases.

\mk

\noindent\emph{Case} 1: $ar\le5040$. If $ar\in\cR$, then since $r\in\cR_2$, relations \eqref{EQ:R2} and \eqref{EQ:R2'} imply $G(ar) \le G(r)$. On the other hand, if $ar\not\in\cR$, then $G(ar) < e^\gamma \le G(r)$. Thus $G(r) \ge G(ar)$ whenever $ar\le5040$.

\mk

Before considering Case~2, we recall that
in \cite[p. 204 (c)]{robin} Robin proved that if $C$ is the largest CA
number with $P(C) < 20000$, then there is no
counterexample $\leq C$ to his inequality \eqref{EQ: robin}. From the property \eqref{k} of
CA numbers, we have $\log C \geq \th(20000),$ where $\th(x)$ is Chebychev's
function.

We also recall that in \cite[p.~359, Corollary~2]{Sch}, Schoenfeld proved that
$$\th(x) > x-\frac{x}{8\log x} \qquad(x\ge19421).$$
A calculation then gives the inequalities
$$\th(20000) > 20000-\frac{20000}{8\log20000} > 19747$$
which, together with Robin's result on $C$, yield the implication
\begin{equation}\label{18990}
5040 < n < e^{19747} \implies G(n) < e^\ga.
\end{equation}

\mk

\noindent\emph{Case} 2: $ar>5040$. If $\log ar < 19747$, then \eqref{18990} gives $G(ar) < e^\ga\leq G(r)$.
On the other hand, if $\log ar \geq 19747$, then from \eqref{EQ: bound} we get
$$G(ar) < e^\ga+\frac{0.6483}{(\log 19747)^2} = 1.787\dotso < 1.790\dotso = \min_{r'\in
  \cR_2} G(r') \le G(r).$$
Thus $G(r) \ge G(ar)$ whenever $ar>5040$.

\mk

This shows that, in both Cases~1 and~2, all elements $r$ of $\cR_2$ are GA2 numbers, so that $\cR_2 \subset \cR'$. Finally, since we already have $\cR_2 \supset \cR'$, we get $\cR_2 = \cR'$. This proves~(i).

\mk 

\noindent(ii). If RH holds, then by Robin's theorem there is no number $n > 5040$ with $G(n) \geq e^\ga$, while from \eqref{GA2>ega} a GA2
number $N$ must satisfy $G(N) \geq e^\ga$.

\mk

\noindent(iii). Let us assume that RH fails. Set
$$\Th:=\sup_{\zeta(\rho)=0} \Re (\rho)$$
so that
$$ 1/2 < \Th \leq 1.$$

Let $N$ denote a CA number of parameter $\vep,$ and define $x=x(\vep)$ by \eqref{defx}.
If $p:=P(N)$ and if $p^+$ is the prime
following $p$, then from \eqref{N=} we have
$$p\leq x_1=x\leq p^+,$$
which implies $x \sim p$ as $N\to \iy$. Further, from
\eqref{PNeqlogN}, we get $p\sim \log N$, which implies
$$x \sim \log N \quad(N\to \iy).$$
In \cite[p. 241]{robinsem}, it is proved that as $N\to \iy$
$$G(N)=e^\ga\left( 1 +\Om_+\left({x^{-b}}\right)\right) \quad (1-\Th
< b <  1/2)$$
which implies that
$$G(N)=e^\ga\left( 1 +\Om_+\left({(\log N)^{-b}}\right)\right) \quad (1-\Th
< b <  1/2).$$
(Here the notation ``$f(N)=\Om_+(g(N))$ as $N\to \iy$'' means that $f(N)>g(N)$ infinitely often, and should not be confused with the notation $\Om(n)$ in \S\ref{SEC: notation}.)
Therefore, there exist infinitely many CA numbers $N$ satisfying $G(N) >
e^\ga$, and, for all~$t$, we have $\max_{n\geq t} G(n) > e^\ga$.

Now we construct two sequences $A_1,A_2,\dotsc$ and $A_1',A_2',\dotsc,$ 
as follows. Let $A_1$ (resp., $A_1'$) be the
smallest (resp., largest)
\footnote{It is highly probable that $A_1=A_1'$. A difficult question is whether $G$ is injective.}
integer $>5040$ such that $G(A_1)=G(A_1')=\mu$. 

Given $i\geq 2,$ assume that $A_1,A_2, \dotsc,
A_{i-1}$ and $A_1',A_2',$ $\dotsc,A_{i-1}'$  have been defined. 
Set $\mu_i:=\max_{n > A_{i-1}'} G(n)$
and let $A_i$ (resp., $A_i'$) be the smallest (resp., largest)
integer $>A_{i-1}'$ with $G(A_i)=G(A_i')=\mu_i$.
Since we have $\mu_i > e^\ga =\limsup G(n)$, infinitely many $A_i$'s can be found.
The numbers $A_i$ are such that 
$$n > A_i \implies G(n) \leq G(A_i)$$
and, therefore, are GA2. 

In the same way, $A$ is proved to be GA2, using $A> 5040$ and $G(A)=\mu$.
To show that $A$ is CA, we apply Lemma \ref{lemNR} with $N_0=55040$, 
$\vep_0=0.03$, $N=A$, and $\vep=1/(\log A \log \log A)$;
since $A$ is GA2 and $A> 5040$, from \eqref{18990} and \eqref{GA2>ega} we
obtain that $N=A> e^{19747} > N_0$. For $n \geq N_0$,
from the definition of $A$ we have $G(n) \leq G(A)$. 
Since $e < N_0 < A$ holds, it follows from Lemma \ref{lemg} that, on the interval $[N_0,+\iy)$, the function $g(t)$ (defined by \eqref{g}) attains its minimum at $t=A$. Thus, for $n\geq N_0$, we have
$$ \frac{\si(n)}{n^{1+\vep}} =G(n) \frac{\log \log n}{n^\vep}=G(n) e^{-g(n)} 
\leq G(A) e^{-g(n)}
\leq G(A) e^{-g(A)}
=\frac{\si(A)}{A^{1+\vep}}$$
and so \eqref{lemNR1} holds. 
Applying Lemma \ref{lemNR} completes the proof of (iii).
\end{dem}

Here is a corollary of the proof of Theorem~\ref{thmGA2}.

\begin{coro} \label{CORlowerbound}
There is no GA2 or extraordinary number between $5041$ and $10^{8576}$.
\end{coro}

\begin{dem}
Since $10^{8576}<e^{19747}$, this follows from \eqref{18990}.
\end{dem}

\section{Comparison between CA and GA1 numbers}\label{secCAGA1}

In this section, we study GA1 numbers. We begin by comparing them with CA numbers.

\subsection{CA and GA1}

By revisiting the proof of \cite[Theorem 3, p. 242]{robinsem}, we shall
prove the following results.
\begin{lem}\label{lemGA11}
Let $N$ be a CA number of parameter $\vep>0$ and assume that
$p:=P(N)\geq 5$. If
\begin{equation}\label{eps1}
\vep > \frac{1}{\log (N/p)\log \log (N/p)},
\end{equation}
then $N$ is also a GA1 number.
\end{lem}

\begin{dem}
Let $q$ be a prime factor of $N$. It follows from \eqref{k} that $6p$
divides $N$ and that $N/q \geq N/p\geq 6 > e$, which implies
$\log \log (N/q) > \log \log e= 0$.
Since $N$ is a CA number, from \eqref{defCA} one has
\begin{equation}\label{eps2}
\frac{\si(N/q)}{(N/q)^{1+\vep} }\leq \frac{\si(N)}{N^{1+\vep}},
\end{equation}
so that
\begin{equation}\label{eps2b}
\frac{\si(N/q)}{\si(N)} \leq \frac{1}{q^{1+\vep}}.
\end{equation}
Since $\log\log N$ and $\log \log (N/q)$ are positive, it follows that
\begin{equation}\label{eps3}
\frac{G(N/q)}{G(N)} \leq \frac{\log\log N}{q^\vep \log\log(N/q)}=
\frac{(N/q)^\vep \log\log N}{N^\vep \log\log(N/q)}=\exp(g(N/q)-g(N)),
\end{equation}
where $g(t)$ is defined by \eqref{g}.
By Lemma \ref{lemg}, using \eqref{epst0} to define $t_0 > e$, we have that $g(t)$ is increasing for $t > t_0$. Now from \eqref{eps1} we deduce that
$$e < t_0 < \frac Np \leq \frac Nq < N$$
and from \eqref{eps3} we get $G(N/q) < G(N).$
This shows that $N$ is GA1.
\end{dem}

\begin{theorem}\label{thm2}
Infinitely many CA numbers are GA1.
\end{theorem}

\begin{dem}
Choose a sufficiently large prime $p$ satisfying
\eqref{P2}, and set $\vep:=F(p,1)$
(so that $x=p$, by \eqref{defx}). Let $N$ be the largest CA number of parameter
$\vep$ (so that $p$ divides $N$, by \eqref{p||N}). From Lemma \ref{lemCA} part (ii) and
\eqref{P2}, we get
$$\log N \geq \th(x)=\th(p) > p + c \sqrt p\log\log\log p,$$
so that
$$\log (N/p)  > p + c \sqrt p\log\log\log p-\log p > p+1.$$
Using the lower bound $\log (1+t) \geq t/(1+t)$, we get
\begin{eqnarray*}
\vep&=&F(p,1) = \frac{\log \left(1+\frac1p\right)}{\log p}\geq \frac{1}{(p+1)\log p}\\
&>& \frac{1}{(p+1)\log (p+1)}>  \frac{1}{\log(N/p)\log\log(N/p)}
\end{eqnarray*}
and Lemma \ref{lemGA11} implies $N$ is GA1.  Since, by Lemma \ref{lemP},
there are infinitely many primes $p$ satisfying \eqref{P2}, the theorem is proved.
\end{dem}

\subsection{CA and not GA1}

To study CA numbers that are not GA1, we need a lemma.

\begin{lem}\label{lemGA12}
Given a prime $p\ge3$, let $N$ be the largest  CA number of parameter $\vep:=F(p,1)$. If
\begin{equation}\label{eps21}
\vep < \frac{1}{\log N\log \log N}\,,
\end{equation}
then $N$ is not  GA1.
\end{lem}

\begin{dem}
As $\vep=F(p,1)$, we have $p^\vep=(p+1)/p=\si(p)/p$. Hence, by \eqref{p||N},
inequality \eqref{eps2} becomes an equality when $q=p$, and so do inequalities
\eqref{eps2b} and \eqref{eps3}. Therefore, with $g$ and $t_0$ defined by \eqref{g} and \eqref{epst0} as in the proof of Lemma \ref{lemGA11}, we get that
$$\frac{G(N/p)}{G(N)} =\exp(g(N/p)-g(N))$$
and, from Lemma \ref{lemg}, that $g(t)$ is decreasing for $t < t_0$. Then \eqref{eps21} implies $N/p < N < t_0$, so
that $G(N) < G(N/p)$. Thus $N$ is not GA1.
\end{dem}

The CA numbers $N$ such that
$P(N) \in \{ 2,3,5,7,11,13,29,59,149\}$ are not GA1. There are two CA
numbers such that $P(N)=23$; the larger one is not GA1, while the
smaller one is GA1. All other CA numbers satisfying $P(N) < 300$ are GA1. (These statements follow by computing all CA numbers $N$ with $P(N)<300$, and calculating those that are GA1---see \S7.)

\begin{theorem}
Infinitely many CA numbers are not GA1.
\end{theorem}

\begin{dem}
Choose a sufficiently large prime $p$ satisfying
\eqref{P1}, and set $\vep:=F(p,1)$
(so that, from \eqref{defx}, $x=p$). Let $N$ be the largest CA number of parameter
$\vep$ (so that, from \eqref{p||N}, $p=P(N)$). From Lemma \ref{lemCA} part (i) and
\eqref{P1}, we get
$$\log N \leq \th(p) +c\sqrt p <  p - c \sqrt p\log\log\log p+c\sqrt p < p,$$
and so
\begin{eqnarray*}
\vep = \frac{\log (1+\frac1p)}{\log p}< \frac{1}{p\log p} <  \frac{1}{\log N\log\log N}.
\end{eqnarray*}
Then Lemma \ref{lemGA12} implies $N$ is not GA1. Since there are
infinitely many primes~$p$ satisfying \eqref{P1}, the theorem is
proved.
\end{dem}

\subsection{Odd GA1 numbers} \label{secodd}

We show that there are infinitely many odd GA1 numbers, and we compute the smallest one.

Let us denote by $\cp_0=\{2,3,5,7,11,13,17,\dotsc\}$ the set of all
primes, and by $\cp$ a subset of $\cp_0$. To $\cp$, we attach the set
$$\cn_\cp:=\{n\ge1:p\mid n \implies p\in \cp\}$$
and the function
$$\th_\cp(x):=\sum_{p\in \cp, \; p\leq x} \log p.$$

A number $N\in \cn_\cp$ is said to be \emph{colossally abundant relative to}
$\cp$ (for short, \emph{CA}$_\cp$) if there exists $\vep > 0$ such that 
$$\frac{\si(N)}{N^{1+\vep}} \geq  \frac{\si(n)}{n^{1+\vep}} \qquad (n\in\cn_\cp).$$
If $M=\prod_{p\in \cp_0} p^{\al_p}$ is an ordinary CA number of
parameter $\vep$, then the factor $N=\prod_{p\in \cp} p^{\al_p}$ is CA$_\cp$,
for the same parameter $\vep,$
and all CA$_\cp$ numbers can be obtained in this way.

\begin{theorem}\label{thm3}
There exist infinitely many odd GA1 numbers.
\end{theorem}

\begin{dem}
First, we observe that Lemma \ref{lemGA11} remains valid if we replace CA with
CA$_\cp,$ for any set $\cp$ with at least 2 elements.

We set $\cp=\cp_0\setminus \{2\}$. The proof of Theorem \ref{thm2}
remains essentially valid. We just have to change the lower bound for
$\log N$ to
$$\log N \geq \th_\cp(p) = \th(p) -\log 2$$
and the inequality $\log (N/p) > p+1$ still holds, so that we may
conclude that $N$ is GA1.
\end{dem}

The smallest CA$_{\cp_0 \setminus \{2\}}$ number that is GA1 is
\begin{align*}
\om:=\ &1058462574572984015114271643676625\\
=\ &3^4\! \cdot\!5^3\! \cdot\!7^2 \!\cdot\!11^2 \!\cdot\! 13\! \cdot\! 17\!\cdot\! 19\! \cdot\!23\!\cdot\! 29\!\cdot\! 31\!\cdot \!37\!\cdot\! 41\!\cdot\! 43\!\cdot \!47\!\cdot\! 53\!\cdot\! 59\!\cdot\! 61\!\cdot\! 67\!\cdot\! 71\!\cdot\! 73.
\end{align*}
From our computation (see \S \ref{subcompR}), $\om$ \emph{is also the smallest odd GA1 number}.

\begin{coro}
There exist infinitely many GA1 numbers that are not SA.
\end{coro}

\begin{dem}
This folllows immediately from \eqref{k} and Theorem \ref{thm3}.
\end{dem}

Of course, the  proof of Theorem \ref{thm3} works for any set of primes $\cp$
such that $\cp_0\setminus \cp$ is finite.

\section{Prime factors of GA1 numbers}\label{secpr}

Here we study prime factors of proper GA1 numbers.

\subsection{An upper bound}\label{secup}

We need the following upper bound.

\begin{theorem}\label{thmplogN}
Given a GA1 number $N$ with $\Om(N) \geq 3$, 
let $p$ be a prime factor of $N$. 
Then for any positive integer $r\le v_p(N)$ we have
$$p\leq (r\log N)^{1/r}\leq \log N.$$
\end{theorem}

\begin{dem}
We have $G(N/p) \leq G(N),$ which implies
\begin{equation}\label{Gp1}
\frac{\si(N/p)N}{(N/p)\si(N)} \leq \frac{\log \log (N/p)}{\log \log N}
=\frac{\log (\log N-\log p)}{\log \log N} \cdot
\end{equation}
Note that $\log \log N > \log \log (N/p) \geq \log \log 4 >0$. We also
have
\begin{eqnarray*}
\log(\log N-\log p) & = & 
\log \left( \log N\left( 1-\frac{\log p}{\log N}\right) \right)\notag\\
  & = &\log \log N + \log \left(1-\frac{\log p}{\log N}\right) \notag
\end{eqnarray*}
so that
\begin{equation}\label{Gp2}
\frac{\log (\log N-\log p)}{\log \log N} =
1-\frac{-\log \left( 1-\frac{\log p}{\log N}\right)}{\log \log N} \cdot
\end{equation}
Further, setting $v=v_p(N)$, the left side of \eqref{Gp1} can be written as
\begin{eqnarray}\label{Gp3}
\frac{\si(N/p)N}{(N/p)\si(N)} & = &
p\ \frac{1+p+\dotsb+p^{v-1}}{1+p+\dotsb +p^v}\notag\\
& = & 1-\frac{1}{1+p+\dotsb +p^v} \geq 1-\frac{1}{1+p+\dotsb +p^r} \cdot
\end{eqnarray}
From \eqref{Gp1}, \eqref{Gp2}, and \eqref{Gp3}, one deduces
\begin{equation}\label{prN}
p^r \leq 1+p+\dotsb +p^r \leq 
\frac {\log \log N}{-\log \left( 1-\frac{\log p}{\log N}\right)}
\leq \frac{\log N \log \log N}{\log p}
\end{equation}
which yields
\begin{equation}\label{Gp4}
p^r \log p \leq\log N \log \log N.
\end{equation}

Let us assume, {\it ab absurdum}, that $p > (r\log N)^{1/r}$. Then we would
have
$$p^r\log p > (r\log N) \frac 1r \log(r\log N)=\log N \log (r\log N)
\geq \log N \log \log N$$
contradicting \eqref{Gp4}. Therefore, $p\leq (r\log N)^{1/r}$
holds. Finally, by calculus, $(r\log N)^{1/r}$ is
decreasing for $r\geq 1$ 
(because $\Om(N)\geq 3$ implies $N\geq 8$ and $\log N > 2$)
and the theorem follows.
\end{dem}

\subsection{Study of $\Om(N)$ where $N$ is GA1} \label{secOm(N)}

We show that there are only finitely many proper GA1 numbers $N$ that have a fixed value of $\Om(N).$

\begin{theorem}\label{thmOmega}
If $k\geq 3,$ then
$$\Pi_k:=\#\{N:N\text{ is GA1 and }\Om(N)=k\}<\infty.$$
\end{theorem}

\begin{dem}
For a GA1 number $N$ with $\Om(N)=k>2,$ let us write $N=p_1p_2\dotsb p_k$ with $p_1 \leq
p_2\leq \dotsb \leq p_k$. We have $N\leq p_k^k$, so that $p_k\geq
N^{1/k}$ holds. But Theorem \ref{thmplogN} yields $p_k \leq \log N$,
whence
$$\frac{\log N}{\log \log N}\leq k$$
and $N$ is bounded. Thus $\Pi_k$ is finite.
\end{dem}

Since $\frac{\log 10^{60}}{\log \log 10^{60}}=28.03\dotsc$, a table of
GA1 numbers up to $10^{60}$ (see \S\ref{seccomp}) allows us to calculate $\Pi_k$
for $k\leq 28$.

We have $\Pi_k=0$ if $3\leq k \leq 12,$ and the following table gives $\Pi_k$ when $13\le k\le28$ (see \cite[Sequence A201558]{oeis}).
$$
\begin{array}{r|cccccccccccccccc}
\hline
k=&13&14&15&16&17&18&19&20&21&22&23&24&25&26&27&28\\
\Pi_k=&2&4&2&1&1&2&4&1&2&3&7&7&7&1&4&7\\
\hline
\end{array}
$$

\subsection{The exponent of the largest prime factor}

First, we observe that the function $t\mapsto 2^t/t$ is an increasing bijection of the interval $[2,+\iy)$ to itself. Let us introduce the inverse function $h$ defined for $x\geq 2$ by
\begin{equation}\label{h}
h(x)=t \qquad \Llr\qquad x=\frac{2^t}{t}\cdot
\end{equation}
We shall need the following lemma.

\begin{lem}\label{lem2pasa}
Let $x$ satisfy $x\geq 2$. Then we have
$2 \leq h(x) \leq 3.08 \log x.$
\end{lem}

\begin{dem}
The lower bound results from the definition of $h$. Let us set $t=h(x)$, so that $x=2^t/t$. By noting that $(\log t)/t\leq 1/e$ holds, we get
$$\frac{h(x)}{\log x} =\frac{t}{t\log 2-\log t}= \frac{1}{\log 2-(\log t)/t}
\leq \frac{1}{\log 2-1/e}=3.0743\ldots$$
which proves Lemma \ref{lem2pasa}. 
\end{dem}

\begin{theorem}\label{coroM}
Let $N$ be a GA1 number with $\Om(N) \geq 3$. Set $R=h(\log N)$, so that $2^R/R=\log N$. Then $N$
divides the number $M=M(N)$ defined by
\begin{equation}\label{M}
M:=\prod_{r=1}^{\lfloor R\rfloor} \;\;
\prod_{((r+1) \log N)^{1/(r+1)} < p\leq (r \log N)^{1/r}} p^r \;\;=\;\;
\prod_{r=1}^{\lfloor R\rfloor} \;\;\prod_{p\leq (r \log N)^{1/r}} p.
\end{equation}
\end{theorem}

\begin{dem}
Since the function $r\mapsto (r \log N)^{1/r}$ is decreasing, this follows from Theorem \ref{thmplogN}.
\end{dem}

For example, if $N=\nu=183783600,$ we compute $R=h(\log \nu)=7.072\dotso$ and find that $M=72\nu$. 

Theorem \ref{coroM} allows the computation of proper GA1 numbers---see \S7.2 and \S7.5.

For the exponent $v_p(N)$ of a prime $p$ in the standard factorization of $N$, Theorem \ref{coroM} provides the upper bound $v_p(N)\leq v_p(M)$, which only depends on the size of $N$.

We now study the exponent of the largest prime factor of a GA1 number.

\begin{theorem}\label{thmvp}
Let $N$ be a GA1 number with $\Om(N)\geq 3,$ and let $p=P(N)$ be its largest
prime factor. Then $v_p(N)=1$.
\end{theorem}

\begin{dem}
Suppose on the contrary that $v:=v_p(N)\geq 2$. Then 
Theorem \ref{thmplogN} implies that $N$ divides the number
$$M_v=M_v(N) :=\left(\prod_{p\leq (v \log N)^{1/v}} p\right)^v \;\;
\prod_{r=v+1}^{\lfloor R\rfloor} \;\;\prod_{p\leq (r \log N)^{1/r}} p$$
with $R$ defined by $2^R/R=\log N.$ 
Thus, from the function $r\mapsto (r \log N)^{1/r}$ being decreasing,
\begin{eqnarray*}
\log N &\leq& \log M_v = v\,\th(v \log N)^{1/v})+
\sum_{r=v+1}^{\lfloor R\rfloor}\th((r \log N)^{1/r})\\
&\leq& 2\th((2 \log N)^{1/2}) +R\th ((3\log N)^{1/3}).
\end{eqnarray*}
From Lemmas \ref{lem2pasa} and \ref{lemmajth}, it follows that
\begin{eqnarray*}
\log N &\leq& 1.000028\left(2\sqrt{2 \log N}
+3.08 \log \log N (3 \log N)^{1/3}\right)\\
&\leq& 2.83\sqrt{\log N}
+4.45 \log \log N (\log N)^{1/3}.
\end{eqnarray*}
Therefore, we have
$$\frac{2.83}{\sqrt{\log  N}}+\frac{4.45}{(\log N)^{2/3} \log \log N}
\geq 1$$
which implies $\log N \leq 15.03$, $N \leq N_0:=3336369$ and $R\leq
6.65,$ so that $N$~must divide $M_v(N_0)$ for some $v$ in the range $2\leq v\leq 6$. But the table
$$
\begin{array}{r|c|c|c|c|c}
\hline
v=&2&3&4&5&6\\
(v \log N_0)^{1/v}=&5.48&3.56&2.78&2.37&2.12\\
M_v(N_0)=&43200=2^6 3^3 5^2&1728=2^6 3^3&64=2^6&64&64\\
\hline
\end{array}
$$
shows that if $v\geq 2$, then the number $M_v(N_0)$ divides $M_2(N_0)=43200$, contradicting the easily-checked fact that none of the 84 divisors 
of $43200$ is a proper GA1 number.
(In fact, we will see in \S\ref{subcompR} that
there is no proper GA1 number $< 183 783 600$.) This proves the theorem.
\end{dem}

\subsection{The largest prime factor of a GA1 number}

For a GA1 number, we now study the largest prime factor itself.

\begin{theorem}\label{thmpequiv}
For GA1 numbers $N$ with $\Om(N)\geq 3,$ the largest prime factor satisfies 
\begin{equation*}
P(N)\sim \log N \quad (N\to\infty).
\end{equation*}
\end{theorem}

\begin{dem}
Let $N$ be a GA1 number satisfying $\Om(N)\geq 3$ and let $p:=P(N)$ be
its largest prime factor. From Theorem \ref{thmplogN}, we know that
\begin{equation}\label{pN}
p\leq \log N.
\end{equation}
It remains to get a lower bound for $p$. The proof
resembles that of Theorem~\ref{thmvp}.

Since $N$ divides $M$ given by \eqref{M} and $p=P(N)$, by Lemma \ref{lemmajth} we have 
\begin{eqnarray}\label{majlogNthp}
\log N & \leq & \th(p)+\sum_{r=2}^R \th((r \log N)^{1/r})\notag\\
&\leq & \th(p) + \th((2\log N)^{1/2}) +\co(R (\log N)^{1/3})\notag\\
&\leq & \th(p)  +\co(\sqrt{\log N}).
\end{eqnarray}
From the Prime Number Theorem and from \eqref{pN}, we get
$$\th(p)=p + \co(p\exp(-c\sqrt{\log p}))=
p + \co(\log N\exp(-c\sqrt{\log \log N})).$$
Therefore, \eqref{majlogNthp} becomes
$$\log N \leq p + \co(\log N\exp(-c\sqrt{\log \log N})),$$
which, together with \eqref{pN}, completes the proof of the theorem.
\end{dem}

 \section{Computation of GA1 numbers}\label{seccomp}

In this section we give several versions of an algorithm to compute GA1 numbers.

\subsection{The Gronwall quotient}\label{subcompG}

We begin with a lemma and a definition.

\begin{lem}\label{lemtest}
Let $n$ be a positive integer with $\Om(n)\geq 3$.
Let $q$ and $p$ be prime factors of $n$ satisfying
$q < p$ and $v_q(n) \leq v_p(n)$. Then we have
$$G(n/q) < G(n/p).$$
\end{lem}

\begin{dem}
We have
\begin{eqnarray*}
\frac{\si(n/q)/(n/q)}{\si(n)/n} & = &
q\ \frac{\si(n/q)}{\si(n)}=
\frac{q+\dotsb + q^{v_q(n)}}{1+q+\dotsb +  q^{v_q(n)}}
=1-\frac{1}{1+q+\dotsb +  q^{v_q(n)}} \\
& \leq & 1-\frac{1}{1+p+\dotsb +
  p^{v_p(n)}} = \frac{\si(n/p)/(n/p)}{\si(n)/n} 
\end{eqnarray*}
which implies $\frac{\si(n/q)}{n/q}\leq \frac{\si(n/p)}{n/p}$.

The lemma follows from $\log \log(n/q) > \log \log (n/p)$, 
since $\Om(n)$ $\geq 3$ implies $\log \log (n/p) \geq \log \log 4 > 0$. 
\end{dem}

We define the {\it Gronwall quotient} $Q(n)$ of a composite integer $n$ to be the
number
$$Q(n):=\max_{\substack{p\mid n\\p\text{ prime}}} \frac{G(n/p)}{G(n)}=
\max_{\substack{p\mid n\\p\text{ prime}}}
\frac{p^{v_p(n)+1}-p}{p^{v_p(n)+1}-1}\;
\frac{\log \log n}{\log \log (n/p)}\cdot$$
GA1 numbers $N$ are characterized by $Q(N)\leq 1$. For example, the ``$Q(r)$'' column in Table~\ref{TABLE: S and G(m)} shows that the only GA1 number $r\in \cR$ is $r=4$.

\mk

Let us introduce a subset $\cs(n)$ of the set of the prime
divisors of~$n$. The elements of $\cs(n)$ are defined by induction. The largest
prime factor of~$n$ is the first element $q_1$ of $\cs(n)$. Now
let us assume that $i\geq 2$ and that the elements $q_1,q_2,\dotsc,q_{i-1}\in \cs(n)$ are known.

If, for all primes $p$ that divide $n$ and are smaller than $q_{i-1}$, we have
$v_p(n)\leq v_{q_{i-1}}(n),$ then there are no further elements of
$\cs(n),$ and we get $\cs(n)\!=\!\{q_1,q_2,\dotsc,q_{i-1}\}$.

If not, then the element $q_i\in \cs(n)$ is defined as the largest prime factor of~$n$ that satisfies
$q_i < q_{i-1}$ and $v_{q_i}(n) > v_{q_{i-1}}(n)$.

From Lemma \ref{lemtest}, if $\Om(n)\geq 3$ we get
$$Q(n)=\max_{p\in \cs(n)} \frac{G(n/p)}{G(n)}=
\max_{p\in \cs(n)} \frac{p^{v_p(n)+1}-p}{p^{v_p(n)+1}-1}\;
\frac{\log \log n}{\log \log (n/p)}\cdot$$

\subsection{A first algorithm}\label{subcomp1}

To compute all proper GA1 numbers $N\leq x$ for a given~$x$, we first
calculate $M=M(x),$ defined by
$$M:=\prod_{r=1}^{\lfloor R\rfloor}  \prod_{p\leq (r \log x)^{1/r}} p,$$
with $R$ such that $(R\log x)^{1/R}=2$. Any GA1 number $N\leq x$ with
$\Om(N)\geq 3$ is a divisor of $M$ (see Theorem~\ref{coroM}).

Thus a first version of the algorithm computes all composite divisors $N$ of~$M,$
and for each of them calculates $G(N/p)/G(N)$ for all $p\in \cs(N)$. If for
some $p\in \cs(N)$ we have $G(N/p)/G(N) > 1$, we stop: $N$ is not GA1. If not,
we compute the Gronwall quotient $Q(N)$ (which involves \emph{all} primes $p$ dividing $N$): $N$ is GA1 if and only if $Q(N)\leq 1.$
\footnote{To avoid roundoff errors, we carry out our computation in floating
point arithmetic with $20$ decimal digits and choose a small $\vep$ (typically,
$\vep=10^{-5}$). In the first step, we keep the $N$'s satisfying $Q(N)
\leq 1+\vep$. For these $N$'s, we start the computation again with
$40$ digits.}.

\subsection{A second algorithm}\label{subcomp2}

A more elaborate version of the algorithm tests only a small number of
the divisors of $M$. First, we define 
$$M_1:=\prod_{r=2}^{\lfloor R\rfloor} \prod_{p\leq (r \log x)^{1/r}} p$$
so that $M_2:=M/M_1$ is squarefree. Let us write $M_2=p_1p_2\dotsb p_s$
where $p_1,p_2,\dotsc ,p_s$ are consecutive primes in ascending order.

As a first step, we compute the set $\cd_0$ of all the composite divisors of $M_1$ and
test each of them for GA1 by the method described above.

A divisor of $M$ whose largest prime factor is
$p_i$ is equal to $d\,p_i$, where $d$ is a divisor of $M$ 
whose largest prime factor is $< p_{i}$. Therefore, we construct 
by induction on $i=1,2,\dotsc,s$ the set $\cd_i'$ containing those divisors
of $M$ whose largest prime factor is $p_i,$ and the set $\cd_i$ 
containing the divisors of $M$ whose largest prime factor is $\leq
p_i.$ Then $\cd_i'$ is equal to $p_i\cd_{i-1}$ and $\cd_i=\cd_i' \cup \cd_{i-1}$.
From Theorem \ref{thmplogN}, for $i=1,2,\dotsc,s$, we only have to
test the elements of $\cd_i'$ that are greater than $\exp(p_i)$. 

\subsection{A third algorithm}\label{subcomp3}

 Let us say that a divisor $d\in \cd_i$ (with
$0\leq i < s$) is {\it bad} if, for every $j$
satisfying $i < j \leq s$, all multiples of $d$ belonging to
$\cd_j$ are smaller than $\exp(p_j)$.

The largest multiple of $d$ belonging to $\cd_j$ is $d\,p_{i+1}p_{i+2}\dotsb
p_j,$ so that $d$ is bad if and only if
$$\log d < \de_i := \th(p_i)+\min_{i< j \leq s} (p_j-\th(p_j)).$$
Therefore, we write $\cg_i \subset \cd_i$ for the set obtained
from $\cd_i$ by removing the bad divisors, i.e., those divisors $d$ satisfying 
$d < \De_i:=\exp(\de_i)$.

Furthermore, we construct $\cg_{i+1}'$ and $\cg_{i+1}$ by removing from
$p_{i+1}\cg_{i}$ and $p_{i+1}\cg_{i}\cup \cg_i,$ respectively, those divisors $d$ that satisfy 
$d < \De_{i+1}=\exp(\de_{i+1})$. 

For $i=1,2,\dotsc,s$, it remains to test the elements of
$\cg_{i}$ whose largest prime factor is equal to $p_i$, that is, the
elements of $\cg_i'$.

\subsection{Results}\label{subcompR}

\emph{The smallest proper GA1 number is}
$$\nu= 183 783 600=2^4\cdot3^3 \cdot 5^2 \cdot 7 \cdot11\cdot13 \cdot17.$$
We compute that $M=M(\nu)=8 \cdot19 \cdot \nu$ and we find that there
is no proper GA1 number $N<\nu$.

\mk

Using the third algorithm, we have  computed all  GA1 numbers $N\le10^{60}$ with
$\Om(N)\geq 3$.

\mk

These results as well as the Maple code can be found 
on the web site \url{http://math.univ-lyon1.fr/~nicolas/GAnumbers.html}.

\mk

We hope to present soon a fourth algorithm, more sophisticated, and able to
compute GA1 numbers up to $10^{120}$.

\section{The number of GA1 numbers up to $x$}\label{secconcl}

Let $Q_1(x)$ be \emph{the number of proper GA1 numbers} $N\leq x$. From
\eqref{M} we know that $Q_1(x)$ does not exceed the number  $\tau(M)$
of divisors of
$$M=M(x):=\prod_{r=1}^{\lfloor R\rfloor} \prod_{p\leq (r \log x)^{1/r}} p$$
with $(R\log x)^{1/R}=2$. It is easy to see that $\log M \sim \log x$ as $x\to\infty,$
and from the estimation of the large values of the function $\tau$
(cf. \cite{hw} or \cite{NR}), it follows that
$$Q_1(x) \leq \exp \left( c\ \frac{\log x}{\log\log x}\right)$$
for some positive $c$.
By estimating the number of {\it good} divisors of $M$ (that is,
divisors that are not bad---see \S \ref{subcomp3}), it might be possible to improve the above
estimate.

It seems more difficult to get a lower bound for $Q_1(x)$. We hope to
return to these questions in another article.\\

\noindent{\bf Acknowledgement} 
The authors thank Kieren MacMillan for his help with TeXing the manuscript.

\def\refname{References}


\begin{thebibliography}{99}

\bibitem{ae} Alaoglu, L., Erd\H{o}s, P.: On highly composite and 
similar numbers. Trans. Amer. Math. Soc. {\bf56}, 448--469 (1944)

\bibitem{briggs} Briggs, K.: Abundant numbers and the Riemann
  hypothesis. Experiment. Math. {\bf15}, 251--256 (2006). \url{http://www.expmath.org/expmath/volumes/15/15.2/Briggs.pdf} (2006). Accessed 23 October 2011

\bibitem{CNS} Caveney, G., Nicolas, J.-L., Sondow, J.: Robin's
  theorem, primes, and a new elementary reformulation of the Riemann
  Hypothesis. Integers {\bf11}, A33.
  \url{http://www.integers-ejcnt.org/l33/l33.pdf} (2011). Accessed 23 October 2011

\bibitem{clms} Choie, Y.-J., Lichiardopol, N., Moree, P., Sole, P.: On Robin's criterion for the Riemann Hypothesis. J. Th\'eor. Nombres Bordeaux {\bf19}, 351--366 (2007). \url{http://arxiv.org/abs/math/0604314} (2006). Accessed 23 October 2011

\bibitem{Dus} Dusart, P.: Estimates of some functions over primes 
without R.H. \url{http://arxiv.org/abs/1002.0442v1} (2010). Accessed 23 October 2011

\bibitem{EN} Erd\H{o}s, P., Nicolas, J.-L.: R\'epartition des nombres superabondants. Bull. Soc. Math. Fr. {\bf 103}, 65--90 (1975). \url{http://www.numdam.org/item?id=BSMF_1975__103__65_0} (1975). Accessed 23 October 2011

\bibitem{gronwall} Gronwall, T.~H.: Some asymptotic expressions 
in the theory of numbers. Trans. Amer. Math. Soc. 
{\bf14}, 113--122 (1913)

\bibitem{hw} Hardy, G.~H., Wright, E.~M.: An Introduction 
to the Theory of Numbers. Heath-Brown, D.~R., Silverman, J.~H. (eds.), 6th ed. Oxford University Press, Oxford (2008)

\bibitem{lagarias} Lagarias, J.~C.: An elementary problem equivalent 
to the Riemann hypothesis. Amer. Math. Monthly {\bf109}, 
534--543 (2002)

\bibitem{Lit} Littlewood, J.~E.: Sur la distribution des nombres
  premiers. C. R. Acad. Sci. Paris S\'er. I Math. {\bf158}, 1869--1872 (1914)

\bibitem{NR} Nicolas, J.-L., Robin, G.: Majorations explicites pour
  le nombre de diviseurs de $N$. Canad. Math. Bull. {\bf
    26}, 485--492 (1983)

\bibitem{ramanujan15} Ramanujan, S.: Highly composite numbers. 
Proc. London Math. Soc. S\'erie 2 {\bf 14}, 347--400 (1915).
Also In: Collected Papers, pp. 78--128. Cambridge University Press, Cambridge (1927)

\bibitem{ramanujan97} Ramanujan, S.: Highly composite numbers, 
annotated and with a foreword by J.-L. Nicolas and G. Robin. 
Ramanujan J. {\bf1}, 119--153 (1997)

\bibitem{robin} Robin, G.: Grandes valeurs de la fonction somme des
  diviseurs et hypoth\` ese de Riemann. J.~Math. Pures Appl.
  {\bf63}, 187--213 (1984)

\bibitem{robinsem} Robin, G.: Sur l'ordre maximum de la fonction somme des
  diviseurs. In: S\'eminaire Delange-Pisot-Poitou Paris 1981-82, pp. 233--242.
  Birkh\"auser, Boston (1983)

\bibitem{Sch} Schoenfeld, L.: Sharper bounds for the
  Chebyshev functions $\th(x)$ and $\psi (x)$. II. 
Math. Comput. {\bf 30}, 337--360 (1976)

\bibitem{oeis} Sloane, N.~J.~A.: The On-Line Encyclopedia of Integer 
Sequences. \url{http://oeis.org} (2011). Accessed 10 December 2011

 \end{thebibliography}
\end{document}